# The Operator Algebras Mentor Network:
Impact of Community-Based Mentoring

**The authors are following the convention to write names in alphabetical order by last name

Anna Duwenig
alduwenig@gmail.com
Department of Mathematics, KU Leuven
Celestijnenlaan 200B
3001 Leuven, Belgium

Kari Eifler
keifler.math@gmail.com

Priyanga Ganesan
pganesan@ucsd.edu
Department of Mathematics, University of California San Diego
La Jolla, CA 92093 USA

Lara Ismert
lara.ismert@gmail.com

Viviana Meschitti
viviana.meschitti@unibg.it
Department of Management, University of Bergamo
24127 Bergamo BG, Italy

Sarah Plosker
ploskers@brandonu.ca
Department of Mathematics and Computer Science
Brandon University
Brandon, MB, Canada R7A 6A9

Karen Strung
strung@math.cas.cz
Institute of Mathematics, Czech Academy of Sciences
Žitna 25
115 67 Prague
Czech Republic

# 1. Abstract

This paper aims to determine if membership within the Operator Algebras Mentor Network (OAMN) is beneficial to its members. The OAMN is an international mentoring initiative that offers support in small groups to women and minority genders in the particularly male-dominated field of operator algebras (OA) in mathematics. Expected advantages of membership include raising awareness of the lack of gender diversity in this field, providing advice to mentees by mentors (e.g., pertaining to career or work/life balance), broadening one's network in OA, etc. A questionnaire was sent to OAMN members and a control group of non-members at similar institutions and similar positions to collect their experience with the mentoring initiative and perception of gender dynamics within the OA discipline, together with basic demographics. The initial analysis of the data collected shows that mentoring junior women and other minority genders in the area has a positive effect on mentees' networking ability, self-promotion, and raising awareness of gender issues within OA as a whole.

# 1. Introduction

## 1.1 Gender Inequality in Mentoring

Mentoring broadly indicates the developmental relationship between a more junior professional (referred to as *protégé* or *mentee* in the literature) and a more experienced one (*mentor*). Mentored individuals experience a variety of benefits, such as increased job satisfaction, less stress, and better self-efficacy compared to non-mentored individuals [1]. Women and racial minorities less often receive informal mentoring [2], so for this reason formal mentoring programs are needed to offer underrepresented groups the opportunity to experience these benefits [3-6]. Formal mentoring programs might even bring positive cultural change in the involved institutions [7].

## 1.2 The Operator Algebras Mentor Network (OAMN)

The Operator Algebras Mentor Network (OAMN) was conceived in 2019 to provide support and empowerment for operator algebraists who are underrepresented or in the minority within mathematics due to gender. Embracing inclusivity, the OAMN uses "underrepresented or in the minority within mathematics due to gender" as an umbrella phrase to encompass various gender identities (e.g., cis-women, transgender, nonbinary). Initially formed to provide support and guidance for individuals facing abuse or harassment, the OAMN has evolved to encompass broader goals aimed at fostering diversity, inclusion, and mentorship within the field of operator algebras.

In alignment with its mission, the OAMN endeavors to retain operator algebraists from such backgrounds, provide early-career professionals with a trusted ally in times of adversity and offer career advice and support outside of traditional institutional frameworks [2]. By forging connections between mathematicians and increasing the visibility of underrepresented voices, the OAMN aims to pave the way for a more equitable and vibrant mathematical community. The structure of the OAMN makes it unique in the panorama of formal mentoring initiatives [9].

With a global footprint spanning 26 countries across North America, Europe, Asia, Africa, the Middle East, and Oceania, the OAMN intends to transcend geographical and institutional boundaries. The OAMN operates through a framework comprising three tiers: Senior Mentors, Junior Mentors, and Mentees. *Senior Mentors* are tenured, tenure-track, or otherwise permanent faculty members or persons in relevant industry positions, normally with at least 6 years of post-PhD experience. *Junior Mentors* are postdoctoral fellows, non-tenured faculty, or persons in relevant industry positions, normally with under 6 years of post-PhD experience. *Mentees* are graduate students or recent graduates who are interested in operator algebras and who belong to the *Mentee Target Group*: reflecting the language used by the OAMN, the authors use the terminology Mentee Target Group (MTG) to refer to the group of mathematicians that are underrepresented or in the minority within mathematics due to gender. This term is applied to Junior and Senior Mentors as well as Mentees.

Members are matched into mentoring groups, henceforth referred to as *OAMN groups*, containing at least one member from each tier and at least one Mentor who is in the MTG. Compared to one-to-one mentoring, group mentoring has the advantage of fostering networking opportunities, and of easing the matching process, because it might be challenging to find mentors from a minority group [11, 12]. In OAMN, Mentees are generally not matched with Mentors at their own institution in an effort to offset a potential power imbalance.

The tiered approach taken by the OAMN is intended to ensure a multifaceted support system. Senior Mentors, who typically hold authority within the community, are intended to step in if cases of abuse or harassment arise. In contrast, Junior Mentors may be able to relate more closely with the Mentees due to their similar academic stages, and therefore similar day-to-day experiences, potentially making Mentees feel more comfortable sharing their thoughts with a peer present. Both Senior and Junior Mentors are able to impart guidance and share experiences to those navigating the early stages of their careers in operator algebras. The expectation is that Junior Mentors likewise receive guidance from Senior Mentors; positioning them in a transitional role between Mentee and Senior Mentor.

The OAMN Board of Directors (Board) typically consists of two members from each of the three tiers for a total of 6 members. All Board members belong to the Mentee Target Group. The current Board elects the next year's Board from a list of vetted candidates, aiming for a matrix of representation in geography (given OAMN's global reach), leadership skills, and experiences.

Each year, the OAMN Board conducts a rematching of the OAMN groups in an effort to promote networking and to accommodate Mentees' changing goals and needs. In addition, the OAMN Board takes on various initiatives aimed at bolstering the professional and personal growth of its members and at making the broader OA community a safer space. These activities have included meetups at conferences, career panels, and mentor-focused workshops.

Integral to the structure of the OAMN are its Bylaws which contain a Code of Conduct. They aim to ensure a safe and respectful environment for all members, and include details about bureaucratic dealings such as the elections for the six Board member positions as well as disciplinary action that will be taken against members who violate the Code of Conduct.

In 2022, the OAMN attained 501(c)(3) non-profit status in the USA which underscored its commitment to its mission of promoting diversity and inclusion within mathematics and gave it financial means to expand their work. In addition, it was awarded an Elsevier Mathematical Sciences Scholarship to fund "initiatives promoting diversity and inclusion in mathematics" [3]. More information about the OAMN can be found in the annual reports [8] as well as an article that appeared in the newsletter of the European Women in Mathematics [4].

At the time of writing, the authors are unaware of any other initiative similar to the OAMN. Some country-specific "women in mathematics" organizations offer one-on-one mentoring that is not specific to a particular field within mathematics (e.g., the Association for Women in Mathematics), and there are several "Women in…" conferences/workshops that are specific to various fields within mathematics, but are more research-based (e.g., the series of "Women in Operator Algebras" workshops, which were held in-person in 2018 and 2023, and in a hybrid format in 2021 at the Banff International Research Station in Banff, AB, Canada). However, the OAMN seems to be unique in offering group-based mentoring that is specific to individuals within operator algebras and is conducted on a worldwide scale.

It should be noted that, although the authors have ties to the OAMN, any opinions stated in this article are those of the authors, and there is no attempt to speak on behalf of the current Board of Directors or any other members of the OAMN.

## 2. Methodology

A survey was created and launched at the end of 2023 targeting all past and present OAMN members as well as general Operator Algebraists outside of the OAMN. The OAMN annual reports are publicly available on the Network's website and contain lists of all past and present members that consented to their names being made public. These lists contained a total of 124 names. For 109 of those 124, the authors were able to locate email addresses. The authors then also compiled a thorough list of other mathematicians in operator algebras not in the OAMN, using resources such as the *Operator Algebra Searchable Information Site* [5] and public attendee lists at operator algebra conferences. Particular effort was made to ensure operator algebraists of various career stages were included (e.g., by specifically searching out conferences aimed at young people to find operator algebraists outside the OAMN but at a comparable career stage to the Mentees of the OAMN), and to find operator algebraists at institutions similar (in geographical region, size, etc.) to the institutions of OAMN members.

The final database contained 241 individuals, 109 of which were at some point a member of the OAMN and the remaining 132 individuals were operator algebraists who were affiliated with an academic institution that was represented by at least one OAMN member. These 132 individuals were used as a control group. The database contained individuals across 75 different institutions located in 20 different countries. The survey consisted of both core questions that were asked of all respondents as well as customized questions for each tier. A unique survey link was emailed to the entire database and surveys were completed anonymously online using *Qualtrics* (a software for advanced qualitative and quantitative research). It should be noted that, due to the small number of individuals of underrepresented genders in OA, all but one author of this article are current or previous members of the OAMN and also took part in the survey.

During the compilation process, early drafts of the survey went through numerous iterations, guided by feedback garnered from several OAMN members after a small test run and feedback from the Research Ethics Committee at the institution of one of the authors [institution's information removed for double blind review process].

The majority of survey questions were framed using a Likert scale, that is, they were of the form "To what extent do you agree with the statement [X]?" with possible answers: *strongly agree, somewhat agree, neither agree nor disagree, somewhat disagree, strongly disagree.*

# 3. Survey Results

## 3.1 Respondent Demographics

The survey was open from December 11th, 2023 through the end of April 2024. During that time, the authors collected a total of 68 respondents with the gender breakdown displayed in Table 1. Of these 68 respondents, 28 respondents are current OAMN members, 35 respondents have never been OAMN members, and 5 respondents were formerly OAMN members.

It should be noted that although every effort was made to garner responses from MTG individuals from outside the OAMN as a control group, all MTG respondents are current or former members of OAMN. This is not surprising—any woman or non-binary person who is interested in being mentored or mentoring individuals in the MTG (i.e., persons who are likely to take the time to fill out the survey) is likely to have been recruited into the OAMN.

Table 1**:** Current position type (rows) by membership and identification with MTG (columns)

|  | **OAMN members (non-MTG and MTG combined)** | **OAMN members: MTG** | **OAMN members: non-MTG** | **Control group (non-OAMN members; all non-MTG)[1]** | **Undisclosed** |
|---|---|---|---|---|---|

---

[1] As mentioned elsewhere in this article, the authors did not intend for the control group to consist entirely of people that do not belong to the Mentee Target Group.

| | | | | | |
|---|---|---|---|---|---|
| Permanent research-focused faculty | 9 | 3 | 6 | 16 | 0 |
| Permanent teaching-focused faculty | 3 | 2 | 1 | 0 | 0 |
| Permanent research-teaching-balanced faculty | 4 | 3 | 1 | 5 | 2 |
| Industry professional | 1 | 1 | 0 | 0 | 0 |
| Research-focused postdoc | 7 | 5 | 2 | 3 | 0 |
| Teaching-focused postdoc | 0 | 0 | 0 | 0 | 0 |
| Graduate student | 8 | 8 | 0 | 7 | 0 |
| Other | 3 | 1 | 2 | 2 | 0 |
| **Total** | **33** | **23** | **10** | **33** | **2** |

The 23 survey respondents of the MTG comprise 21 women and two nonbinary people.

## 3.2 Value of Networking

Overwhelmingly, the respondents believe that networking with academics from other institutions is important. Of the 68 respondents,

- 97% agree that networking with academics from other institutions is important,
- 96% agree that a global network is beneficial in pursuing their career goals as an operator algebraist, and
- 85% agree that discussing their career interests in a group of peers will benefit them in the long run.

When asked "How many mathematicians from the MTG have you had communication with?", 85% of OAMN respondents had communicated with fewer than eleven mathematicians from the MTG; to be precise, 38% had communicated with up to 5 MTG mathematicians, and 47% with 5-10. The control group of non-OAMN members had nearly identical statistics (38% communicating with 0-5 MTG mathematicians, and 44% communicating with 5-10 MTG mathematicians). For respondents who have been or currently are members of the OAMN, the authors were unsurprised to find that, most respondents (24 out of 27) strongly agree or somewhat agree with the statement "Having a member from the MTG in my group enables me to broaden my perspective and understanding of academia.".

## 3.3 Insights on Gender Inclusion and Community Belonging

Regardless of gender, 82% of OAMN respondents *strongly agree* with the statement "It is important to raise awareness about instances of gender-based harassment within the operator algebras community", while 71% of respondents in the control group expressed the same level of agreement. The survey shows there is resounding support within the OA community for increasing awareness of issues related to gender equity. However, responses from the MTG to the statement "I feel supported with respect to any perceived gender barriers" indicate that this sentiment does not equate to perceived support: only 14% of MTG respondents *strongly agreed*, while 41% *somewhat agreed*, 36% *neither agreed nor disagreed*, and 9% *somewhat agreed*, with no respondents strongly disagreeing. One respondent, who has never been a member of the OAMN, wrote

> *"I am a senior fully promoted person, also cis white male ... I am sympathetic, but surely naive about EDI issues."*

Responses to this question from individuals within the MTG are similar to their responses to "I 'fit in' within the operator algebras mathematical community" and "I do not feel isolated in the operator algebras community."

In response to the question, "How aware are you about issues that mathematicians from MTG face?", 12% of men responded *very much*, 71% *somewhat*, and 17% *not at all*. Out of the control group respondents, i.e., those who are not and have never been an OAMN member, 74% responded *somewhat*, while current OAMN members were more evenly split between *very much* and *somewhat* (46% and 50%, respectively). This demonstrates that involvement in the OAMN is correlated to one's awareness of issues faced by individuals in the MTG. An OAMN member wrote,

> *"I am a senior mentor, and the only cis man in our group of five ... I have found our meetings enormously beneficial to me in terms of hearing the concerns of junior members of the community."*

The responses to the statement "Being an operator algebraist is part of my identity" varied drastically by gender. While 31% of men strongly agreed, only 4% of MTG strongly agreed.

Conversely, 22% of MTG somewhat disagreed while only 2% of men somewhat disagreed. A full breakdown is shown in Figure X.

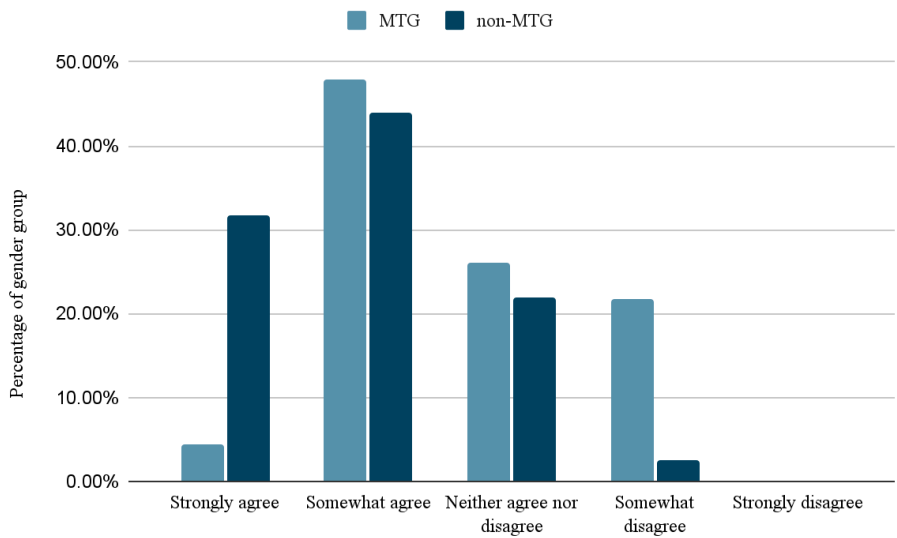

Figure 1: Responses to the question "Being an operator algebraist is part of my identity."

When asked to rate their comfort level in speaking up at a conference, 76% of men feel *highly* or *moderately comfortable,* while only 38% of MTG respondents feel *highly or moderately comfortable*. In contrast, all 21 mentees and junior mentors polled rated a higher comfort level when asked the question "How comfortable do you feel speaking up during an OAMN meeting?" with 70% of respondents feeling *very comfortable* and 30% of respondents feeling *somewhat comfortable*.

## 3.4 Access to Support Networks

Regarding access to a support system or network for MTG at their respective home institutions, 43% of MTG respondents and 36% of non-MTG respondents report they are aware of such a support system. This could lead to a disadvantage for MTG grad students or postdocs—especially those with a non-MTG supervisor—as they may not be made aware of such support systems even if they do exist.

Even without a designated support program/network for MTG at one's home institution, access to role models of the MTG, within or outside of one's institution, has a positive impact on the retention of a member of the MTG [15, 16].

In addition to the OAMN providing mentorship opportunities specifically to early-career MTG mathematicians in operator algebras, the OAMN is also unique in that its mentors are likewise within the operator algebras community. Only 3 of the 23 total graduate student or postdoc respondents reported having a mentor in operator algebras, other than their supervisor, at their home institution. Similarly, only 9% of all respondents reported being aware of a formal mentoring program solely for operator algebraists at their home institution. Access to networking opportunities within one's particular area of research provides unique opportunities for research and career advancement.

## 3.5 The Organization of the OAMN

The authors aim to evaluate the effectiveness of the organizational structure. To this end, additional questions were asked to the 33 OAMN members, consisting of 11 Mentees, 12 Junior Mentors, and 10 Senior Mentors. The survey data reveals both strengths and areas for improvement, supporting some choices made by OAMN's early Board of Directors and guiding future refinements.

Participants rated the value of having "access to a mentor network that is specific to operator algebras" at 8.7 out of 10. This underscores the significance of a tailored approach to mentorship in a niche field where institutional support may be limited. One respondent emphasized that "global mentoring is crucial for operator algebras" due to the limited number of specialists at many institutions.

OAMN groups are typically composed of three to five members. Comments in the survey responses suggested that smaller groups of three members (one Senior Mentor, one Junior Mentor, and one Mentee) are optimal for fostering in-depth discussion, with one Mentor in a group of 5 people expressing the following:

> *"Hearing diverse opinions from people in different positions and from different backgrounds has been very helpful. However, sometimes I am less likely to express something (personal, for example) because I am talking to a bigger group."*.

While the data supports the decision to favor smaller group configurations, average satisfaction was still high among groups with more members:Satisfaction ratings for group size aligned with this trend:

Table 2: Responses to the questions "how satisfied are you with the size of your current mentor group" and "the size of my group is appropriate to enable an in-depth conversation"

| Group Size | Respondents | Avg. Satisfaction of OAMN group size | Avg. "Enabled In-Depth Conversations" |
|---|---|---|---|
| 3 | 3 | 4.7 | 4.5 |
| 4 | 16 | 4.3 | 4.4 |
| 5 | 11 | 4.3 | 4.1 |

Unfortunately, forming such small groups may be challenging due to the limited number of Senior Mentor volunteers compared to the higher number of Mentees and Junior Mentors.

One question that arose in analyzing the data was whether or not a member's meeting frequency affects their satisfaction with OAMN. Members who met with their group once per month reported the highest average satisfaction rating of 4.2. In contrast, those meeting once every two months had an average satisfaction rating of 3.7, which was comparable to the 3.8 rating from members who selected "Less Often." These results suggest that monthly meetings are associated with the highest levels of satisfaction. Based on this, the authors recommend organizing OAMN groups to meet on a monthly basis. One member emphasized the importance of proactive scheduling, stating:

> *"I would prefer if the mentors were more proactive with finding times to meet/chat."*

However, global mentoring is not without its challenges. Scheduling meetings across time zones can be difficult, and Mentors may lack familiarity with the academic systems of their Mentees' institutions. As one participant pointed out, this can lead to *"generic and conservative advice"* in some cases. Despite these limitations, the consensus is that the benefits—such as freedom from local institutional politics and access to a broader network—far outweigh the drawbacks.

The inclusion of members from the MTG in each group enabled participants to broaden their perspective and understanding of academia, with a majority of respondents (55%) *strongly agreeing* and another 33% *somewhat agreeing*. As one member noted:

> *"The diverse cultural background really enriches the conversations and gives us a better idea of how everything works in different stages."*

One of the more profound impacts of OAMN has been its role in addressing gender-based harassment in academia. Participants reported increased confidence in identifying and addressing such issues. One respondent shared:

> *"Before, I felt like pointing out any of these situations was pointless, as it seemed like nobody actually cared about it. I do not feel like that anymore."*

Others highlighted how the network provides moral support and a collective acknowledgment of these challenges, reducing the isolation often felt by individuals.

The global nature of OAMN mentorship was widely praised. Respondents valued the diverse experiences and perspectives offered by members from different institutions and countries. One mentee explained,

> *"It has given me more of a bird's-eye view of academia on a global scale."*

## 3.6 Impacts of the OAMN

### 3.6.1 Impact on Mentees

The OAMN aims to provide tailored support, encouragement, and access to valuable opportunities that can shape mentees academic and professional trajectories. One mentee shared a particularly poignant story:

> *"I received an email from my then OAMN mentor reminding me of the deadline for the PhD position I currently have. Without this encouragement, I would have not applied and I would have likely quit academia without achieving my goal of getting a PhD."*

During rematching, efforts are made to ensure Mentees are paired with at least one Mentor from their continent, recognizing that some topics related to an institution's geographical context can require more localized expertise. When asked whether the Mentors in their group are knowledgeable about advising on topics specific to their geographical location, 45% of Mentees *strongly agreed*, 22% *somewhat agreed*, 22% *neither agreed nor disagreed*, and 11% *somewhat disagreed*.

Many Mentees reported significant academic and professional benefits through their OAMN participation. One Mentee described their Senior Mentor as "the closest mathematician to my research area," noting how personalized feedback on CVs and talks, as well as networking advice, proved invaluable. Another Mentee emphasized the collaborative opportunities enabled by the Network:

> *"Two other members of my first mentor group ended up authoring a paper together, and this has led to a really strong research program. I also got advice and ample encouragement to apply for grants."*

Survey results further reflect this positive experience, with 80% of Mentees rating their participation in the OAMN as 4 or 5 out of 5.

The OAMN also serves as a conduit for information on job openings, speaking engagements, and funding opportunities. Responses varied widely, with some respondents reporting exposure to anywhere from zero to twelve job opportunities, zero to eight speaking opportunities, and zero to six funding opportunities. This disparity highlights an inconsistency in the provision of such information, as the Board of the OAMN does not make blanket announcements about opportunities. Instead, Mentors play a key role in sharing relevant information within their groups. This variability suggests that while some Mentors actively provide Mentees with resources, others may not.

### 3.6.2 Impact on Mentors

While the primary mission of the OAMN is to support early-career operator algebraists and retain them in the field, the Network also aims to enhance mentorship skills among mathematicians and increase the visibility of underrepresented or minority operator algebraists. The survey responses from Mentors highlight the significant personal and professional growth they experience through their involvement in the program.

The authors were very glad to see that 100% of OAMN Mentor respondents (i.e., both Junior and Senior Mentors) said they had made changes to how they organize conferences, advise students, or choose collaborators as a result of their participation in the OAMN. This indicates that the program is successfully cultivating more reflective and effective mentors. One Junior Mentor noted:

> *"It has made me more confident about my abilities, and I have often been encouraged by my mentor groups to apply for certain academic opportunities."*

This reflection underscores how mentorship can be a mutually beneficial process, offering Mentors opportunities for self-assessment and growth while positively influencing their Mentees.

For some Mentors, the OAMN has provided a sense of community and encouragement during difficult times in their own academic journeys. A Junior Mentor shared their experience:

> *"At a time when I was not doing well due to being in a challenging temporary academic position, I felt supported."*

# 4. Conclusion, Future Work, Opinions

## 4.1 Conclusions

The authors designed and conducted a survey studying various factors involving women and minority genders in the field of operator algebras—a field in which these groups are disproportionately underrepresented, in comparison to the broader area of mathematics, in which they are already underrepresented (see Section 3). The survey aimed to collect experiences with the mentoring initiative and perceptions of gender dynamics within the operator algebra discipline. The final database contained 68 individuals, with 28 being current or former OAMN members and 40 serving as a control group (see Section 4.1).

Research has shown that mentored individuals experience various benefits such as increased job satisfaction, less stress, and better self-efficacy compared to non-mentored individuals (see Section 2.1). The OAMN was established in 2019 to support and empower operator algebraists who are underrepresented or in the minority within mathematics due to gender. The OAMN operates globally, transcending geographical and institutional boundaries, and uses a tiered mentoring approach to provide a multi-faceted support system. The Network spans 26 countries and includes Senior Mentors, Junior Mentors, and Mentees (see Section 2.2).

The survey results show that 97% of respondents agree that networking with academics from other institutions is important and 96% agree that a global network is beneficial in pursuing their career goals as an operator algebraist. Additionally, 85% agree that discussing their career interests in a group of peers will benefit them in the long run (see Section 4.2.). The survey reveals strong support for raising awareness about gender-based harassment within the OA community, with 77% of OAMN respondents and 71% of control group respondents strongly agreeing. However, only 14% of MTG respondents strongly agree that they feel supported with respect to perceived gender barriers (see Section 4.3.).

Regarding access to support networks, 43% of MTG respondents and 36% of non-MTG respondents report being aware of such support systems at their home institutions. The OAMN provides unique mentorship opportunities specifically for early-career MTG mathematicians in operator algebras (see Section 4.4.).

The survey also included additional customized questions for the OAMN members. The survey data reveals that participants rate the value of having access to a mentor network specific to operator algebras at 8.7 out of 10. Smaller groups of three members (one Senior Mentor, one Junior Mentor, and one Mentee) are considered optimal for fostering in-depth discussions (see Section 4.5.).

The OAMN provides tailored support that significantly impacts mentees' academic and professional trajectories, with 80% of Mentees rating their participation as 4 or 5 out of 5. Success stories include mentees receiving critical encouragement to pursue PhD positions, invaluable networking advice, collaborative opportunities leading to published research, and access to personalized feedback—though the variability in sharing job and funding opportunities suggests room for improvement in mentor engagement (see Section 4.6.1). The OAMN not only supports early-career operator algebraists but also fosters significant professional growth among mentors, with 100% of Mentors reporting changes in how they organize conferences, advise students, or choose collaborators. This highlights OAMN's success in cultivating reflective and effective mentorship (see Section 4.6.2).

## 4.2 Future Work

Although the authors spent considerable time and effort in designing the survey and crafting the survey questions, some questions could be developed further. In particular, one question asks respondents to rate the degree to which they agree with the statement "I feel supported with respect to any perceived gender barriers." It would be of interest to ask a follow-up question, asking how often they feel supported, to get a more nuanced answer. Additionally, OAMN member respondents felt that male Senior Mentors were very beneficial while male Junior Mentors were only somewhat beneficial; that is, there is a clear preference for male Senior Mentors over male Junior Mentors.

There is potential to conduct such a survey regularly (perhaps biannual). The survey results are meant to help guide the direction and operation of the OAMN in a meaningful, data-driven manner. Follow-up interviews could further explore members' responses regarding OAMN experience. For example, while Mentors feel that OAMN is only a "3" in its effectiveness, whereas Mentees were much more satisfied, on the whole, with their experience in OAMN. This warrants further investigation so as to optimize the satisfaction and effectiveness of the OAMN for its Mentors.

The authors were surprised to find that MTG responses to the question "How supported do you feel with respect to issues related to gender equality in OA?" had a distinctive mode of "Agree" rather than "Strongly Agree," given that all MTG survey respondents have been or are currently members of the OAMN. Upon reflection, the authors recognize that even being part of OAMN does not unilaterally resolve systemic gender barriers within the mathematical community, at one's home institution, or within academia. Further questions related to this topic, including an examination of possible correlation between responses to this question and respondents' overall satisfaction with OAMN, could be done.

The authors plan to also investigate how the length of one's membership in the OAMN affects their satisfaction with the OAMN. Survey data does not indicate a correlation when considering all members of OAMN, but the authors want to consider a possible correlation by subgroup within the OAMN (Mentees, Junior Mentors, Senior Mentors).

## 4.3 Selected Comments from the Survey

Below are selected suggestions we received and felt were worth sharing for other research areas who may be looking to begin an organization like the OAMN:

1. I would only comment that the OAMN is a hugely valuable and powerful resource for young people and very worthwhile; but at the same time, **the organisers and board need to be mindful of their own careers and emotional wellbeing and mental health when determining how much time and energy they can commit to this role**. It is possible for something to be very much worth doing and very valuable, and at the same time for you personally to have given as much as you can to it and need to step back from it… **One thing that the OAMN does not (to my knowledge) explicitly provide is direct support and mentoring for the role of board member.**
2. **I am a senior mentor, and the only cis man in our group of five.** I hope I have been useful to the junior mentors (I think it helps that I am senior enough to have sat on grant panels, organized quite a few conferences, have been a department chair etc.). More selfishly, **I have found our meetings enormously beneficial to me in terms of hearing the concerns of junior members of the community.** The perspective of members who are at different types of institution is also useful to me.
3. I am very happy with all the advice given so far, I think that all this aid could have been more beneficial more halfway of my phd, than now, I think I lose many opportunities due to my ignorance and cant [sic] fix those now. With this, I would suggest to promote this group to more phds as it is very helpful and most of us never had a mentor before.

4. OAMN is such a great resource. **I wish my peers in other areas of math had such an opportunity.** I feel so grateful that this organization exists!
5. I was assigned a mentor group in the middle of last year and have yet to have a meeting with them. **I would prefer if the mentors were more proactive with finding times to meet/chat.**
6. It would be nice to figure out a way for **more in person meetups!**
7. Many of my colleagues outside the OAMN often want to look to the OAMN as a way to reach the mentee group in order to advertise conferences, workshops, or opportunities or to look for speakers. **It would be nice if there were an easier way to connect mentees with these opportunities outside the OAMN.**
8. Some sort of group-based research workshop would be great, if possible to organize
9. One can be or feel excluded without belonging to any identified minority group. In particular, it can be the case that within the community, the perception exists that "this person will be fine", whereas this is not the case at all. Diversity is often invisible.
10. **I am a senior fully promoted person, also cis white male, and answered accordingly. I now mentor more than need mentoring. I am sympathetic, but surely naive about EDI issues.**

Select comments to the question, "How has being a member of the OAMN affected your willingness or comfort in pointing out instances of gender-based harassment?" are below with key takeaways in bold font:

1. I feel more comfortable and more confident that reports will be taken seriously.
2. …I think/hope I would speak up, and **this would possibly be easier knowing I would if needed have support from OAMN**
3. Before I felt like pointing out any of these situation was pointless, as it seem like nobody actually care about it . I do not feel like that anymore and I feel that something can be done in a more formal way.
4. It is useful to me to be able to point to things we have discussed in the OAMN. **The mere existence of the OAMN also provides evidence that gender-based harassment occurs,** and that others (not just me!) acknowledge it. In particular, I don't feel like a lone voice crying in the wilderness.

Select responses to the prompt, "List the benefits and challenges of **global** mentoring (in contrast to institutional-based mentoring)" are included below, with key takeaways in bold font:

1. …The diverse approaches and ideas given helped the discussion as well as their experiences, **even to consider or not moving to other countries based of their experiences**…
2. I think having global mentoring is crucial for operator algebras. Since this is a relatively small field, **many people do not have many other operator algebraists in their home institution**. I also feel like I have grown personally from talking to people from different countries and cultures. It has given me more of a bird's eye view of academia on a global scale…

3. Get **valuable connections and insight into the academic job market in other countries**, as it is often very different.
4. The main benefit is that **people feel completely confident that they can discuss issues and problems they are having at their home institution without fear that this will impact how they are treated at that institution.** The main challenge (beyond timetabling/organisation) is that **the mentors in the group may have little or no understanding of either the academic system or the interpersonal relationships at the relevant institution** - this can lead them to feel constrained to very generic and conservative/careful advice, and also make them feel relatively powerless to help, beyond giving moral support, in situations where assistance/support is needed.
5. It has been very useful to me as **I am in a small institution, which while friendly, is not big enough to have support groups in my specific area.** It has helped me make many contacts, friends and even some potential collaborators.
6. I identify more strongly with the international operator algebra community than I do with the broader mathematics community… I also greatly appreciate that I am able to be in contact with people from countries like India, where there is a huge operator algebra community, but who don't always travel to European and North American conferences that I am more likely to go to.

# Acknowledgements

The authors would like to thank the Brandon University Research Ethics Committee (BUREC) for their time and advice. Qualtrics license was purchased via support from the Canada Research Chairs (CRC) Program grant number 231250 (S.P.).